# A FEW RESULTS ON THE INFIMUM OF REGULAR POLYGONS EQUAL-SIZE SPLIT LINE

Yuyang Zhu
Department of Math and Physics, Hefei University 230601
email: zhuyy@hfuu.edu.cn

**Abstract**: If an n-side unit regular polygon is divided into m equal sized parts, then what is the minimum length of the split line $l_{m,n}$? This problem has its practical application in real world. This paper proved that $l_{2,3} = \sqrt{\dfrac{\sqrt{3}\pi}{12}}$, $l_{3,3} = \dfrac{\sqrt{3}}{2}$, and $\dfrac{1}{2}\sqrt{n\pi \operatorname{ctan}\dfrac{\pi}{n}} \le \lim\limits_{m\to\infty}\dfrac{l_{m,n}}{\sqrt{m}} \le \sqrt{\dfrac{\sqrt{3}}{2}n\operatorname{ctan}\dfrac{\pi}{n}}$, with constructor methods in functional analysis and combinatorial analysis.

**MR(2000) Subject Classification  52C35**

## 1. Introduction

The study on discrete combinatorial geometry has been heated in recent years due to its application in programming and game theory and the development of computer science [1]. The following problem will be investigated in this paper: A square house is subdivided into several equal-area parts. What is the best strategy of subdivision to minimize the subdividing material used? As the dividing lines get shorter, the use of subdividing materials is less. In addition, in the building of ship hulls, lengths of welding lines should be minimized to save welding material and improve hulls strength. Therefore investigation on the infimum of subdividing lines' length is very important.

Let $A_n$ be an n-vertices unit regular polygon, with an area of $\dfrac{n}{4}c\tan\dfrac{\pi}{n}$. If $A_n$ is divided into m equal-sized parts, then let $L_{m,n}$ be the length of split lines in $\operatorname{int} A_n$ (interior of $A_n$), and its minimum $l_{m,n} = \inf L_{m,n}$. Variational methods in functional analysis are used to prove that $l_{2,3} = \sqrt{\tfrac{\sqrt{3}\pi}{12}}$, $l_{3,3} = \tfrac{\sqrt{3}}{2}$, and constructor methods are used to prove that:

$$l_{4,3} \le \sqrt{\dfrac{3\sqrt{3}\pi}{8}}, \quad l_{4,4} \le 2, \quad l_{6,3} \le \dfrac{1}{2}\sqrt{\sqrt{3}\pi} + \dfrac{3}{2}\left(\sqrt{3}-\sqrt{\dfrac{\sqrt{3}}{\pi}}\right),$$

$$l_{n+1,n} \le n\sqrt{\dfrac{1}{n+1}} + \dfrac{n}{2}\left(\cos\dfrac{\pi}{n} - \sqrt{\dfrac{1}{n+1}}\right)\csc\dfrac{\pi}{n}, \quad l_{m,n} \ge \dfrac{1}{2}\left(\sqrt{mn\pi c\tan\dfrac{\pi}{n}} - n\right).$$

## 2. Lemma and Proofs

**Lemma 2.1(see[2])** *Let $X, Y$ be both normed space, while $Y^*$ be the dual space of $Y$. Also define $R$ to be real number space with target function $f: X \to R$, and constraint function $g: X \to Y$. For*

$$\min_{x \in X} f(x), \quad s.\ t.\ g(x) = 0. \tag{2.1}$$

*If $x_0$ is a local optimized solution of (2.1), $f(x)$ is differentiable at $x_0$, $g(x)$ is a $C^1$ mapping at*

*the vicinity of $x_0$, and $g'(x_0)$ is a surjection, then $\exists \lambda \in Y^*$, such that*

$$f'(x_0) + \lambda g'(x_0) = 0. \tag{2.2}$$

**Lemma 2.2** *Assume $M(t)$ is continuous on $J = [a,b]$, if arbitrary $\rho(t)$ which has continuous second derivative on $J = [a,b]$, equals zero at $a, b$ and $\int_a^b M(t)\rho(t)dt = 0$, Then $M(t) \equiv 0$ on $J = [a,b]$.*

**Proof** Assume $M(t) \neq 0$ at the vicinity of $t_0$. Let $M(t) > 0$, $\xi_1 < t < \xi_2$. If

$$\rho_0(t) = \begin{cases} (t-\xi_1)^4(t-\xi_2)^4, & \xi_1 \leq t \leq \xi_2, \\ 0, & t \notin [\xi_1, \xi_2]. \end{cases}$$

Since $\rho_0''(t)$ is continuous on $J = [a,b]$ and $\rho_0(t) = 0$ at $a,b$, and according to the previous assumption,

$$0 = \int_a^b M(t)\rho_0(t)dt = \int_{\xi_1}^{\xi_2} M(t)\rho_0(t)dt.$$

However, the integral on the right hand side is positive everywhere which makes its result positive, and is inconsistent with previous derivation. The lemma is thus proved. □

**Lemma 2.3** *Let functional $f(x) = \int_a^b F(t, x, x')dt$, with field of definition $D(J)$ composed of functions of the following properties:*

(i) $x$ has continuous second derivative.

(ii) $x(a) = \alpha, x(b) = \beta$, and $g(x) = \int_a^b G(t, x, x')dt = k$ ($k$ is a constant). $F$ and $G$ both have continuous partial derivatives. If $f(x)$ has its extreme value at $x_0(t)$, and $G_x' - \dfrac{d}{dt}G_{x'}'$ is continuous but not constant zero, then $\exists \lambda$, such that $x_0(t)$ is a solution to Euler equation

$$F_x' + \lambda G_x' - \frac{d}{dt}(F_{x'}' + \lambda G_{x'}') = 0. \tag{2.3}$$

**Proof** Let $\rho_1(t)$ and $\rho_2(t)$ be functions with continuous second derivative and

$$\rho_i(a) = \rho_i(b) = 0 \ (i=1,2).$$

Consider the function

$$g(\varepsilon_1, \varepsilon_2) = \int_a^b G(t, x_0 + \varepsilon_1\rho_1 + \varepsilon_2\rho_2, x_0' + \varepsilon_1\rho_1' + \varepsilon_2\rho_2')dt - k,$$

Take derivative of on the kernel function. Since $\int_a^b [\rho_2 G_{x'}']'dt = [\rho_2 G_{x'}']\big|_a^b = 0$, and using integration by parts, we have

$$g_{\varepsilon_2}'(\varepsilon_1, \varepsilon_2) = \int_a^b (G_x'\rho_2 + G_{x'}'\rho_2')dt = \int_a^b (G_x' - \frac{d}{dt}G_{x'}')\rho_2 dt.$$

Taking into account that $G_x' - \dfrac{d}{dt}G_{x'}'$ is not constant zero, according to the proof of lemma 2.2, we can choose $\rho_2(t)$ so that $g_{\varepsilon_2}'(0,0) \neq 0$. Consider another function

$$f(\varepsilon_1,\varepsilon_2) = \int_a^b F(t, x_0 + \varepsilon_1\rho_1 + \varepsilon_2\rho_2, x_0' + \varepsilon_1\rho_1' + \varepsilon_2\rho_2')dt.$$

According to our assumption, it has extreme value at $(0,0)$ given the condition $g(\varepsilon_1,\varepsilon_2)=0$.
According to lemma 2.1, there must exist a constant $\lambda$ such that

$$\frac{\partial f}{\partial \varepsilon_i} + \lambda \frac{\partial g}{\partial \varepsilon_i} = 0.$$

Now that $\rho_2(t)$ has been chosen, thus we consider the case of $i=1$ only. It is obvious that

$$\begin{aligned}\frac{\partial f}{\partial \varepsilon_1} + \lambda \frac{\partial g}{\partial \varepsilon_1} &= \int_a^b (F_x'\rho_1 + F_{x'}'\rho_1')dt + \lambda \int_a^b (G_x'\rho_1 + G_{x'}'\rho_1')dt \\ &= \int_a^b (F_x' - \frac{d}{dt}F_{x'}')\rho_1 dt + \lambda \int_a^b (G_x' - \frac{d}{dt}G_{x'}')\rho_1 dt \\ &= \int_a^b [(F_x' + \lambda G_x') - \frac{d}{dt}(F_{x'}' + \lambda G_{x'}')]\rho_1 dt = 0.\end{aligned}$$

Where $\rho_1(t)$ can be any function with continuous second derivative while equals zero at $a,b$.
According to lemma 2.2

$$F_x' + \lambda G_x' - \frac{d}{dt}(F_{x'}' + \lambda G_{x'}') = 0.$$

The lemma is proved. □

**Lemma 2.4** *Let* $L \in D(L)$ *($D(L)$ is the set of all smooth curves) be a smooth curve, and* $L_1$ *be a line segment with fixed length. If the area encircled by $L_1$ and $L$ is fixed, then the length of $L$ is minimized when it is an arc.*

**Proof** As shown in Figure 1, the area encircled by fixed length line segment $AB$ and smooth curve L is $\int_a^b x\,dt$, and the length of L is $\int_a^b \sqrt{1+x'^2}\,dt$.

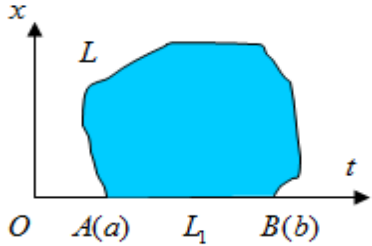 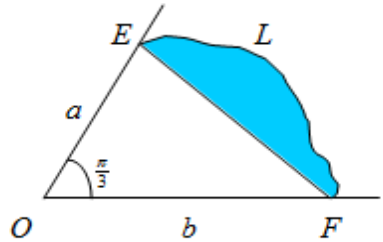

**Fig. 1** Diagram showing proof of lemma 4    **Fig. 2** Diagram showing proof of lemma 5

Let $g(t,x,x') = \int_a^b x\,dt - k$, $f(t,x,x') = \int_a^b \sqrt{1+x'^2}\,dt$. Since $\int_a^b x\,dt = k$ which is a constant, $F(t,x,x') = \sqrt{1+x'^2}$, $G(t,x,x') = x$. According to (2.3) in lemma 2.3,

$$0 + \lambda - \frac{d}{dt}\left(\frac{x'}{\sqrt{1+x'^2}}\right) = 0, \text{ thus}$$

$$x = \int \frac{\lambda t + C_1}{\sqrt{1-(\lambda t + C_1)^2}} dt$$. Therefore $(x-C_2)^2 + \left(t + \frac{C_1}{\lambda}\right)^2 = \frac{1}{\lambda^2}$, where $C_1$ and $C_2$ are both constants.

In other words, the curve with minimized length must be an arc. The lemma is proved. □

**Inference 2.4.1**(Duality of isoperimetric) *A close domain with fixed area has its boundary length minimized when the boundary is a circle.*

**Lemma 2.5** *Consider a sector with angle $\frac{\pi}{3}$, a vertex $O$ and a moving curve connecting two moving points $E, F$ at its corners. If the area of the sector is $A$, which is a constant, then if and only if $L$ is an arc with radius $OE$, will the length of $L$ be minimized, which is $\sqrt{\frac{2}{3} A \pi}$.*

**Proof** As shown in Figure 2, let $OE = a$, $OF = b$, then the area of $\triangle EOF = \frac{\sqrt{3}}{4} ab$. The area encircled by $L$ and $EF$ is then $A - \frac{\sqrt{3}}{4} ab$. According to lemma 2.4, if and only if $L$ is an arc will its length be minimized. Let its radius be $r$, center be $O'$, $D$ be the middle point of $EF$, $|EF| = 2d$, and $\angle DO'E = \angle DO'F = \alpha$, as shown in Figure 3. The area $\sigma$ of the segment and the length $s$ of arc $L$ can be calculated. According to law of cosines, $2d = \sqrt{a^2 + b^2 - ab}$, hence

$$s = 2r \arcsin \frac{\sqrt{a^2 + b^2 - ab}}{2r},$$

$$\sigma = r^2 \arcsin \frac{\sqrt{a^2 + b^2 - ab}}{2r} - \frac{1}{4} \sqrt{a^2 + b^2 - ab} \sqrt{4r^2 - a^2 - b^2 + ab}.$$

Let $f(a,b,r) = 2r \arcsin \frac{\sqrt{a^2 + b^2 - ab}}{2r}$,

$$g(a,b,r) = r^2 \arcsin \frac{\sqrt{a^2 + b^2 - ab}}{2r} - \frac{1}{4} \sqrt{a^2 + b^2 - ab} \sqrt{4r^2 - a^2 - b^2 + ab} - A + \frac{\sqrt{3}}{4} ab.$$

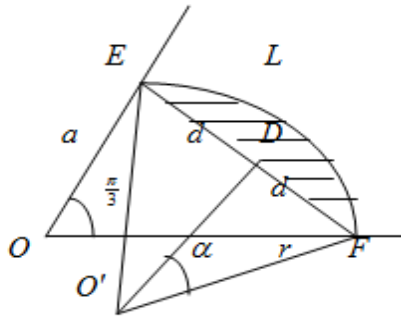
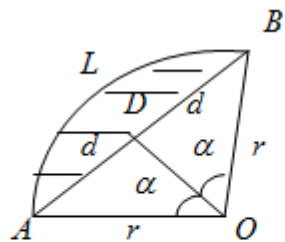

**Fig. 3 Diagram showing proof of lemma 5**  **Fig.4 Diagram showing proof of lemma 6**

then $g(a,b,r) = 0$. According to the conditional extreme value theorem, consider the following equations

$$\begin{cases} \dfrac{\partial f}{\partial r} + \lambda \dfrac{\partial g}{\partial r} = 0, \\ \dfrac{\partial f}{\partial a} + \lambda \dfrac{\partial g}{\partial a} = 0, \\ \dfrac{\partial f}{\partial b} + \lambda \dfrac{\partial g}{\partial b} = 0, \\ g = 0. \end{cases}$$

The first equation gives $1 + \lambda r = 0$. Substitute that in the second and third equation and we get $r = a = b$. Substitute this in $\sin\alpha = \dfrac{\sqrt{a^2 + b^2 - ab}}{2r}$ gives $2\alpha = \dfrac{\pi}{3}$, which makes the center of the arc at $O$. Furthermore, $g(a,b,r) = 0$, gives $a = r = \sqrt{\dfrac{6}{\pi}A}$, then $s = 2\alpha r = \sqrt{\dfrac{2}{3}A\pi}$. On the other hand, when $a \to 0$, due to the fixed area of the section, $b$ can approach infinity, and the same with $d = \dfrac{1}{2}\sqrt{a^2 + b^2 - ab}$ and the length of $L$. Therefore the extreme value point obtained must be the minimized point. The Lemma is proved. □

**Lemma 2.6** *For circular segments with fixed area $A$, the arc length can be minimized, if the string $AB$ with length $2d$ and the radius of the segment meets relation $r = d$. The minimum value of the arc length is $\sqrt{2A\pi}$.*

**Proof** Consider the case when central angle is equal or less than $\pi$. As shown in Figure 4, let the arc be $L$, string be $AB$ with length $2d$ and middle point D, $OD \perp AB$, $\angle AOD = \angle BOD = \alpha$, radius be $r$, and the area $\sigma$ and arc length $s$ be respectively

$$\sigma = r^2 \arcsin\dfrac{d}{r} - d\sqrt{r^2 - d^2} = A, \quad s = 2r\arcsin\dfrac{d}{r}.$$

Let $f(r,d) = 2r\arcsin\dfrac{d}{r}$, $g(r,d) = r^2\arcsin\dfrac{d}{r} - d\sqrt{r^2 - d^2} - A = 0$. According to the conditional extreme value theorem, consider the following equations:

$$\begin{cases} \dfrac{\partial f}{\partial r} + \lambda \dfrac{\partial g}{\partial r} = 0, \\ g = 0, \\ \dfrac{\partial f}{\partial d} + \lambda \dfrac{\partial g}{\partial d} = 0. \end{cases}$$

The first equation gives $1 + \lambda r = 0$, and by substituting that in the third equation gives $r = d$. Substituting that in the second equation gives $r = d = \sqrt{\dfrac{2A}{\pi}}$, $s = 2r\arcsin\dfrac{d}{r} = 2\sqrt{\dfrac{2A}{\pi}}\cdot\dfrac{\pi}{2} = \sqrt{2A\pi}$.

When the central angle is not less than $\pi$, we can prove in the similar way that when $d = r$, $L$ reaches its minimum $\sqrt{2A\pi}$. The lemma is proved. □

**Lemma 2.7(see[3])** *A regular simplex $\sum_{(n+1)}$ in $n$ dimensional Euclidean space has edge length $\rho$, and $M$ is an arbitrary point in the simplex. Then the sum of distances from $M$ to each side of the simplex $d_i$ is $\sum_{i=0}^{n} d_i = \left(\dfrac{n+1}{2n}\right)^{\frac{1}{2}} \rho$.*

**Inference 2.7.1** *The sum of distances from an arbitrary point in a regular triangle to its sides is* $\frac{\sqrt{3}}{2}$.

## 3. The side length infimum of equiareal division line for unit regular triangles

**Theorem 3.1** $l_{2,3} = \sqrt{\frac{\sqrt{3}\pi}{12}}$.

**Proof** (i) If the equiareal division line has zero or one intersection with the unit regular triangle, then the split line must be a closed curve. The enclosed area is $\frac{\sqrt{3}}{8}$, or one half of the area of $A_3$. Let $L$ be this closed curve, with length $s$. According to lemma 2.4, $s$ is minimized if $L$ is an arc. Assume the radius of the arc is $r$, then $\pi r^2 = \frac{\sqrt{3}\pi}{8}$, $s \geq 2\pi r = \sqrt{\frac{\sqrt{3}\pi}{2}} > \sqrt{\frac{\sqrt{3}\pi}{12}}$, namely $L_{2,3} > \sqrt{\frac{\sqrt{3}\pi}{12}}$.

(ii) If $L$ has two intersections with one side of $A_3$, then mark them as $A, B$. According to lemma 2.4, $L$ must be an arc when its length $s$ is minimized. Furthermore, according to lemma 2.6, $s \geq \sqrt{2A\pi} = \sqrt{\frac{\sqrt{3}\pi}{4}} > \sqrt{\frac{\sqrt{3}\pi}{12}}$, namely $L_{2,3} > \sqrt{\frac{\sqrt{3}\pi}{12}}$.

(iii) If $L$ has two intersections with two sides of $A_3$, then again mark them as $A, B$. Mark the common vertex of these sides as $O$, then the area encircled by $\angle AOB$ and $L$ has an area of $\frac{\sqrt{3}}{8}$ (half of the area of $A_3$). According to lemma 2.5, $s \geq \sqrt{\frac{2}{3}\frac{\sqrt{3}}{8}\pi} = \sqrt{\frac{\sqrt{3}\pi}{12}}$, hence $L_{2,3} \geq \sqrt{\frac{\sqrt{3}\pi}{12}}$. All possible cases of equiareal divisionting have been considered in the cases above. Therefore $l_{2,3} = \inf L_{2,3} \geq \sqrt{\frac{\sqrt{3}\pi}{12}}$. On the other hand, according to (iii) and lemma 2.5, when $OA = OB = \sqrt{\frac{3\sqrt{3}}{4\pi}}$ and $\angle AOB = \frac{\pi}{3}$, the arc length $s = \sqrt{\frac{\sqrt{3}\pi}{12}}$, and the area of circular sector $AOB$ is $\frac{\sqrt{3}}{8}$ ($A, B$ are two intersection points on two sides of $A_3$). In this case $L$ is an arc connecting $AB$ with length $s = \sqrt{\frac{\sqrt{3}\pi}{12}}$. thus $l_{2,3} = \inf L_{2,3} \leq \sqrt{\frac{\sqrt{3}\pi}{12}}$. Concluding that, $l_{2,3} = \inf L_{2,3} = \sqrt{\frac{\sqrt{3}\pi}{12}}$. The theorem is proved. □

**Theorem 3.2** $l_{3,3} = \frac{\sqrt{3}}{2}$.

**Proof** (i) If equiareal division line $L$ has none or only one intersection with regular unit triangle, then mark part of $L$ ($L_1$) which encloses $\frac{1}{3}$ area of $A_3$, namely $\frac{\sqrt{3}}{12}$. A circle with area $\frac{\sqrt{3}}{12}$ has perimeter of $\sqrt{\frac{\sqrt{3}\pi}{3}} > \frac{\sqrt{3}}{2}$. According to lemma 2.4 and inference 2.4.1, the length of $L_1$ is $s_1 \geq \sqrt{\frac{\sqrt{3}\pi}{3}} > \frac{\sqrt{3}}{2}$. Thus $L_{3,3} > \frac{\sqrt{3}}{2}$.

(ii) If $L$ has two intersections with sides of $A_3$, then as shown in figure 5(a) and 5(b). According to the proof of (i), $s > L_1$, namely $s_1 > \frac{\sqrt{3}}{2}$. Thus $L_{3,3} > \frac{\sqrt{3}}{2}$.

(iii) If $L$ has two intersections with one side of $A_3$, but one intersection with another side, as shown in Figure 6. $L = L_1 \cup L_2 \cup L_3$. Because $L$ evenly split $A_3$ by three, the area enclosed by $L_1, L_2$ and $AB$ is $\frac{1}{3}$ of the area of $A_3$, namely $\frac{\sqrt{3}}{12}$. According to lemma 2.4 and lemma 2.6, the length of $L >$ the length of $L_1 + L_2 \geq \sqrt{2\frac{\sqrt{3}}{12}\pi} > \frac{\sqrt{3}}{2}$. Namely $L_{3,3} > \frac{\sqrt{3}}{2}$.

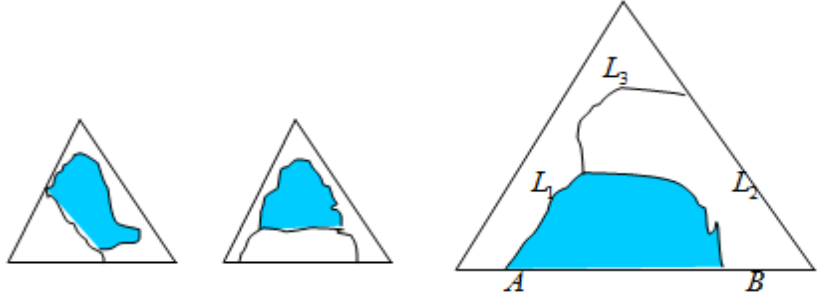

Fig. 5(a)     Fig. 5(b)     Fig. 6

(iv) If $L$ has three intersections with one side of $A_3$ and no intersection with other sides, as shown in Figure 7, then the area enclosed by $L_1$, $L_2$ and $AB$ is $\frac{1}{3}$ of the area of $A_3$, namely $\frac{\sqrt{3}}{12}$. Similar to the proof of (iii), according to lemma 2.4 and 2.6, the length of $L >$ the length of

$L_1 + L_2 \geq \sqrt{2\frac{\sqrt{3}}{12}\pi} > \frac{\sqrt{3}}{2}$.

Hence $L_{3,3} > \frac{\sqrt{3}}{2}$.

(v) If $L$ has one intersection with each side of $A_3$, as shown in Figure 8, then
$L = L_1 \cup L_2 \cup L_3$, $L_1 \cap L_2 \cap L_3 = D$.

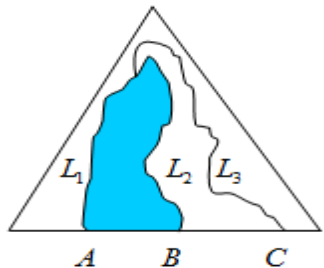 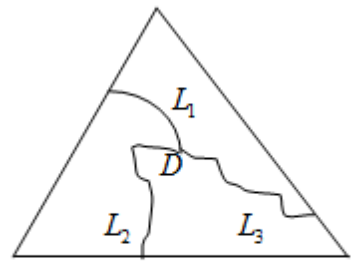

Fig. 7     Fig. 8

According to inference 2.7.1, since $A_3$ is a regular simplex, the distance sum from $D$ to all sides of $A_3$ is $\frac{\sqrt{3}}{2}$. Hence the length of $L \geq \frac{\sqrt{3}}{2}$, and $L_{3,3} \geq \frac{\sqrt{3}}{2}$.

(vi) If $L$ has 4 intersections with $A_3$, as shown in Figure 9(a) or Figure 9(b).

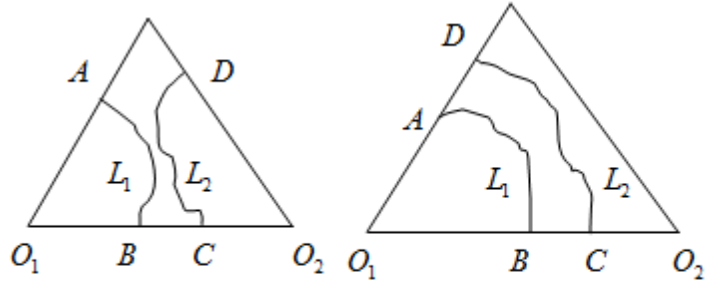

Fig. 9(a)    Fig. 9(b)

The area enclosed by $\angle AO_1B$ and $L_1$, along with $\angle CO_2D$ and $L_2$ are both $\frac{1}{3}$ the area of $A_3$. According to lemma 2.5, the lengths of $L_1$ and $L_2$ are $\geq \sqrt{\frac{2}{3} \cdot \frac{\sqrt{3}}{12}\pi} = \sqrt{\frac{\sqrt{3}}{18}\pi}$, and the length of $L >$ the length of $L_1 + L_2 \geq \sqrt{2\frac{\sqrt{3}}{12}\pi} > \frac{\sqrt{3}}{2}$. In the case of 9(b), according to lemma 2.5, since the length of $L_2 >$ the minimum length of $L_1$, the minimum length of $\geq \sqrt{\frac{\sqrt{3}}{18}\pi}$. Thus the length of $L >$ the length of $L_1 + L_2 \geq \sqrt{2\frac{\sqrt{3}}{12}\pi} > \frac{\sqrt{3}}{2}$, and $L_{3,3} > \frac{\sqrt{3}}{2}$.

As all possible cases have been considered above $l_{3,3} = \inf L_{3,3} \geq \frac{\sqrt{3}}{2}$.

On the other side, if $L$ is coincidentally the connection from the center of $A_3$ to middle points of sides of $A_3$, then $L$ evenly split $A_3$ by three, and $L$ has a length of $\frac{\sqrt{3}}{2}$, hence $l_{3,3} = \inf L_{3,3} \leq \frac{\sqrt{3}}{2}$.

To conclude, $l_{3,3} = \inf L_{3,3} = \frac{\sqrt{3}}{2}$. The theorem is proved. □

**4. An asymptotic estimation of the minimum and maximum length of equiareal division line**

**Theorem 4.1** $l_{m,n} \geq \frac{1}{2}\left(\sqrt{mn\pi\operatorname{ctan}\frac{\pi}{n}} - n\right)$.

**Proof** By doing an equiareal division of a unit n-vertices polygon $A_n$ into $m$ parts. Each division has an area of $\frac{n}{4m}\operatorname{ctan}\frac{\pi}{n}$. According to the inference 2.4.1 of lemma 2.4, the boundary length is minimized when the boundary of each division is a circle. If the area of a circle is $\frac{n}{4m}\operatorname{ctan}\frac{\pi}{n}$, its perimeter is $\sqrt{\frac{n\pi}{m}\operatorname{ctan}\frac{\pi}{n}}$. The total perimeter of $m$ circles is $\sqrt{mn\pi\operatorname{ctan}\frac{\pi}{n}}$.

Subtracted by the side length of $A_n(n)$, we get $\sqrt{mn\pi\operatorname{ctan}\frac{\pi}{n}}-n$. On the other hand, the equiareal division boundaries must be intersections of neighboring boundaries, thus the total length should be halved, namely $L_{m,n} \geq \frac{1}{2}\left(\sqrt{mn\pi\operatorname{ctan}\frac{\pi}{n}}-n\right)$.

Hence $l_{m,n} = \inf L_{m,n} \geq \frac{1}{2}\left(\sqrt{mn\pi\operatorname{ctan}\frac{\pi}{n}}-n\right)$. The theorem is proved. □

According to theorem 3, $\lim\limits_{m\to\infty}\frac{l_{m,n}}{\sqrt{m}} \geq \frac{1}{2}\sqrt{n\pi\operatorname{ctan}\frac{\pi}{n}}$.

Since equal sized regular hexagons can occupy the whole plane, we can put $m-t$ $(m,t \in N, t \leq m)$ equal sized regular hexagons $B_1, B_2, \cdots, B_{m-t}$ inside regular polygon $A_n$, and meet the condition $B = \bigcup_{k=1}^{m-t} B_k \subset \operatorname{int} A_n$ (where $\operatorname{int} A_n$ represents the internal of $A_n$), $\partial B \cap \partial A_n = \varnothing$, and the area of $B_k$ $S(B_k) = \frac{n}{4m}\operatorname{ctan}\frac{\pi}{n}$ $(k=1,2,\cdots,m-t)$. Since $B$ is a single communication, $\partial B$ is a polygon, and sides of the polygon have equal length (the same as the side length of $B_1$). Choose an arbitrary side on $\partial B$ and draw a regular hexgon $B_{m-t+1}$ of which area is equal to that of $B_1$. If that always meets the condition $\partial B_{m-t+1} \cap \partial A_n \neq \varnothing$, then the maximum embedded number of regular hexagons is $m-t$ for $m$ equiareal division of $A_n$. Given that condition is met, equiareally divide ring shaped $A_n - B$ with t line segments represented by $E_1F_1, E_2F_2, \cdots, E_tF_t$, where $E_k \in \partial A_n$, $F_k \in \partial B$ $(k=1,2,\cdots,t)$. Denote the total length of these line segments by $l_t$, then we have the following theorem.

**Theorem 4.2** *If the maximum embedded number is $m-t$ for $m$ equiareal division of $A_n$, then*

$$\lim_{m\to\infty}\frac{l_t}{m} = 0.$$

**Proof** Since $m-t$ is the maximum embedded number for $m$ equiareal division of $A_n$, by definition, $\forall E \in \partial A_n$, $\inf\limits_{F\in\partial B}\{|EF|\} \leq$ diameter of regular hexagon $= \sqrt{\frac{2n}{3\sqrt{3}m}\operatorname{ctan}\frac{\pi}{n}}$. According to the definition of $l_t$, $l_t \leq t \cdot$ diameter of $B_1 = t \cdot \sqrt{\frac{2n}{3\sqrt{3}m}\operatorname{ctan}\frac{\pi}{n}}$. On the other hand, $0 \leq t \leq m$, thus

$$0 \leq \lim_{m\to\infty}\frac{l_t}{m} \leq \lim_{m\to\infty}\frac{t}{m}\cdot\sqrt{\frac{2n}{3\sqrt{3}m}\operatorname{ctan}\frac{\pi}{n}} = 0.$$

The theorem is proved. □

An asymptotic estimation can be made on the supremum of $l_{m,n}$.

**Theorem 4.3** $\lim\limits_{m\to\infty} \dfrac{l_{m,n}}{\sqrt{m}} \leq \sqrt{\dfrac{\sqrt{3}}{2}n\mathrm{ctan}\dfrac{\pi}{n}}$ .

**Proof** $B$ and $l_t$ have been defined previously. The ring shape $A_n - B$'s area

$$S(A_n - B) \leq \text{perimeter of } A_n \text{ times diameter of } B_1 = n \cdot \sqrt{\dfrac{2n}{3\sqrt{3}m}\mathrm{ctan}\dfrac{\pi}{n}}, \text{ thus}$$

$$\lim\limits_{m\to\infty} S(A_n - B) \leq \lim\limits_{m\to\infty} n \cdot \sqrt{\dfrac{2n}{3\sqrt{3}m}\mathrm{ctan}\dfrac{\pi}{n}} = 0.$$

If $m$ is sufficiently great, there exists a positive real number $a$, such that $l_t < a \cdot$ perimeter of $A_n$. Hence, there exists $N_0 \in N^+$, for all $m > N_0$, $l_t < a \cdot n$. We can conclude that if $m > N_0$, then

$$l_{m,n} \leq l_t + \dfrac{1}{2}(m-t) \cdot \text{perimeter of } B_1 + \dfrac{1}{2}\text{perimeter of } B$$

$$\leq l_t + \dfrac{1}{2}m \cdot \text{perimeter of } B_1 < a \cdot n + \dfrac{1}{2} \cdot m \cdot 6\sqrt{\dfrac{n}{6\sqrt{3}m}\mathrm{ctan}\dfrac{\pi}{n}},$$

thus $\lim\limits_{m\to\infty} \dfrac{l_{m,n}}{\sqrt{m}} \leq \sqrt{\dfrac{\sqrt{3}}{2}n\mathrm{ctan}\dfrac{\pi}{n}}$. The theorem is proved. □

**Note 1** According to theorem 3 and 5, $\dfrac{1}{2}\sqrt{n\pi\mathrm{ctan}\dfrac{\pi}{n}} \leq \lim\limits_{m\to\infty} \dfrac{l_{m,n}}{\sqrt{m}} \leq \sqrt{\dfrac{\sqrt{3}}{2}n\mathrm{ctan}\dfrac{\pi}{n}}$. Since $\left|\dfrac{\sqrt{3}}{2} - \dfrac{\pi}{4}\right| < 0.1057$, the estimation has considerable accuracy.

The length of curve connecting two points on a plane is minimized if the curve is a straight line, and the length infimum of equiareal division line for a regular polygon is related to the occupation or packing of equal size graphs. For a given area $n$-vertices regular polygon, as $n$ increases its perimeter decreases. Only regular triangle, rectangle and hexagons can fully occupy a plane, so for a large $m$, to minimize the division lines, a $n$-vertices regular polygon should be divided by regular polygons as much as possible.
There are two guesses following:

**Guess 1** $\lim\limits_{m\to\infty} \dfrac{l_{m,n}}{\sqrt{m}} = \sqrt{\dfrac{\sqrt{3}}{2}n\mathrm{ctan}\dfrac{\pi}{n}}$ .

In some special cases,

**Guess 2** $l_{4,3} = \sqrt{\dfrac{3\sqrt{3}\pi}{8}}$ ; $l_{6,3} = \dfrac{1}{2}\sqrt{\sqrt{3}\pi} + \sqrt{3} - \dfrac{1}{2} \cdot \sqrt{\dfrac{3\sqrt{3}}{\pi}}$ ; $l_{3,4} = 1 + \dfrac{2}{3}$ ; $l_{4,4} = 2$ ;

$l_{n+1,n} = n\sqrt{\dfrac{1}{n+1}} + \dfrac{n}{2}\left(\cos\dfrac{\pi}{n} - \sqrt{\dfrac{1}{n+1}}\right)\csc\dfrac{\pi}{n}$ $(n > 3)$.